\documentclass[letterpaper, 10 pt, conference]{ieeeconf}  

\usepackage{graphicx}  
\usepackage{epstopdf}  
\usepackage{subfig}    
\usepackage{amsmath}   
\usepackage{boxhandler}
\usepackage{dsfont}    
\usepackage[flushleft]{threeparttable} 
\let\labelindent\relax 
\usepackage{enumitem}  
\usepackage{booktabs}  
\usepackage{multirow}  
\usepackage{url} 	   
\usepackage{setspace}  
\usepackage{textcomp}  
\usepackage{color}     
\usepackage{booktabs} 
\usepackage{array}    
\usepackage{caption}  
\usepackage[font=footnotesize]{caption}
\usepackage{epsfig}   
\usepackage{mathptm}  
\usepackage{amssymb} 	
\usepackage{amsxtra}    
\usepackage{verbatim}  
\usepackage{graphics} 
\usepackage{pdflscape} 
\usepackage{enumerate}
\usepackage{enumitem}
\usepackage{setspace}
\usepackage[export]{adjustbox}
\usepackage{setspace}
\usepackage{bm}
\usepackage{lipsum}
\usepackage{algorithm}
\usepackage[noend]{algpseudocode}
\usepackage[noadjust]{cite} 

\IEEEoverridecommandlockouts                              

\newcommand*{\TitleFont}{%
       \usefont{\encodingdefault}{\rmdefault}{m}{n}%
       \fontsize{24}{29}%
       \selectfont}
       
\begin{document}      
     

\title{\TitleFont
Experimental Modeling of Cyclists Fatigue and Recovery Dynamics Enabling Optimal Pacing in A Time Trial
}

\author{Faraz Ashtiani \quad Vijay Sarthy M Sreedhara \quad Ardalan Vahidi \quad Randolph Hutchison \quad Gregory Mocko
%
\thanks{Faraz Ashtiani ({\tt\footnotesize fashtia@g.clemson.edu}), Vijay Sarthy M Sreedhara ({\tt\footnotesize vsreedh@g.clemson.edu}), Ardalan Vahidi ({\tt\footnotesize avahidi@clemson.edu}),  and Gregory Mocko ({\tt\footnotesize gmocko@clemson.edu}) are with the Department of Mechanical Engineering, Clemson University, Clemson, SC 29634-0921, USA. Randolph Hutchison ({\tt\footnotesize randolph.hutchison@Furman.edu}) is with the Department of Health Science, Furman University, Greenville, SC, 29613-1000, USA.}%
}

\maketitle
\thispagestyle{empty}
\pagestyle{empty}

\begin{abstract}
Improving a cyclist performance during a time-trial effort has been a challenge for sport scientists for several decades. There has been a lot of work on understanding the physiological concepts behind it. The concepts of Critical Power ($CP$) and Anaerobic Work Capacity ($AWC$) have been discussed often in recent cycling performance related articles. $CP$ is a power that can be maintained by a cyclist for a long time; meaning pedaling at or below this limit, theoretically, can be continued for infinite amount of time. However, there is a limited source of energy for generating power above $CP$. This limited energy source is $AWC$. After burning energy from this tank, a cyclist can recover some by pedaling below $CP$. In this paper we utilize the concepts of $CP$ and $AWC$ to mathematically model muscle fatigue and recovery of a cyclist. Then, the models are used to formulate an optimal control problem for a time trial effort on a 10.3 km course located in Greenville SC. The course is simulated in a laboratory environment using a CompuTrainer. At the end, the optimal simulation results are compared to the performance of one subject on CompuTrainer.
\end{abstract}


\section{Introduction}
\label{introduction}
Fatigue due to prolonged exercise can be defined as a decline in muscle performance, which accompanies a sensation of tiredness \cite{abbiss2005models,bigland1984changes}. Thus, fatigue causes an inability to exercise at the required intensity. Exercise intensity falls under one of severe, high, or moderate categories. This classification is based on either change rate of blood lactate levels \cite{rose2007quantitative}, maximum oxygen uptake ($\dot{V}O_{2,max}$)\cite{hall2002affective,welch2007affective}, and power output \cite{welch2007affective}. Critical power remarks severe and heavy intensity domains \cite{black2016muscle} and \cite{poole1988metabolic}. 

The notion of ``critical power" was introduced by the authors in \cite{monod1965work}. They defined it from the notions of \textit{maximum work} and \textit{maximum time of work}, where there is a border power between aerobic and anaerobic exercise. One can define critical power as this boundary at which a human can generate power for an infinite amount of time \cite{monod1965work}. While pedaling below or above $CP$, a cyclist is utilizing aerobic or anaerobic metabolic systems, respectively. There is a limited amount of energy in human body to perform anaerobic exercise. When this energy is finished, cyclists cannot pedal at a power level above $CP$. This limited energy reservoir is called Anaerobic Work Capacity ($AWC$). $AWC$ is similar to the notion of \textit{maximum work} discussed in \cite{monod1965work}. Researchers in some studies suggest an experimental method for determining $CP$ and $AWC$ which takes each subject multiple lab visits \cite{bull2000effect,bergstrom2014differences,gaesser1995estimation,hill1993critical,housh1990methodological}. During these visit days, subjects ride on a certain power level ($P$) above $CP$ until they cannot hold it at the same level anymore, and the time they can hold their power is measured ($\Delta t$). The data from all of the days can be used to calculate $AWC$ as below:

\begin{equation}
AWC = (P - CP)\Delta t
\label{eq_basic-model}
\end{equation}

To reduce the lab visit days, a 3 minute all-out test is developed in \cite{vanhatalo2007determination}. As its name suggests, subjects are required to continuously sprint and exert their maximum power for three minutes. As shown in figure \ref{fig_3mao}, the subject's average power during the last 30 seconds on the test is considered to be $CP$. Also, the area between exerted power and $CP$ in a power vs. time plot is the subject's $AWC$. This means a cyclist's $AWC$ can last for 2 and a half minutes in a maximal effort. The accuracy of this test is validated in \cite{vanhatalo2008robustness} and \cite{johnson2011reliability}.

\begin{figure}
\centering
\setlength\fboxsep{0pt}
\setlength\fboxrule{0pt}
\setlength\abovecaptionskip{-4pt}
\setlength\belowcaptionskip{-20pt}
\fbox{\includegraphics[width= 8cm]{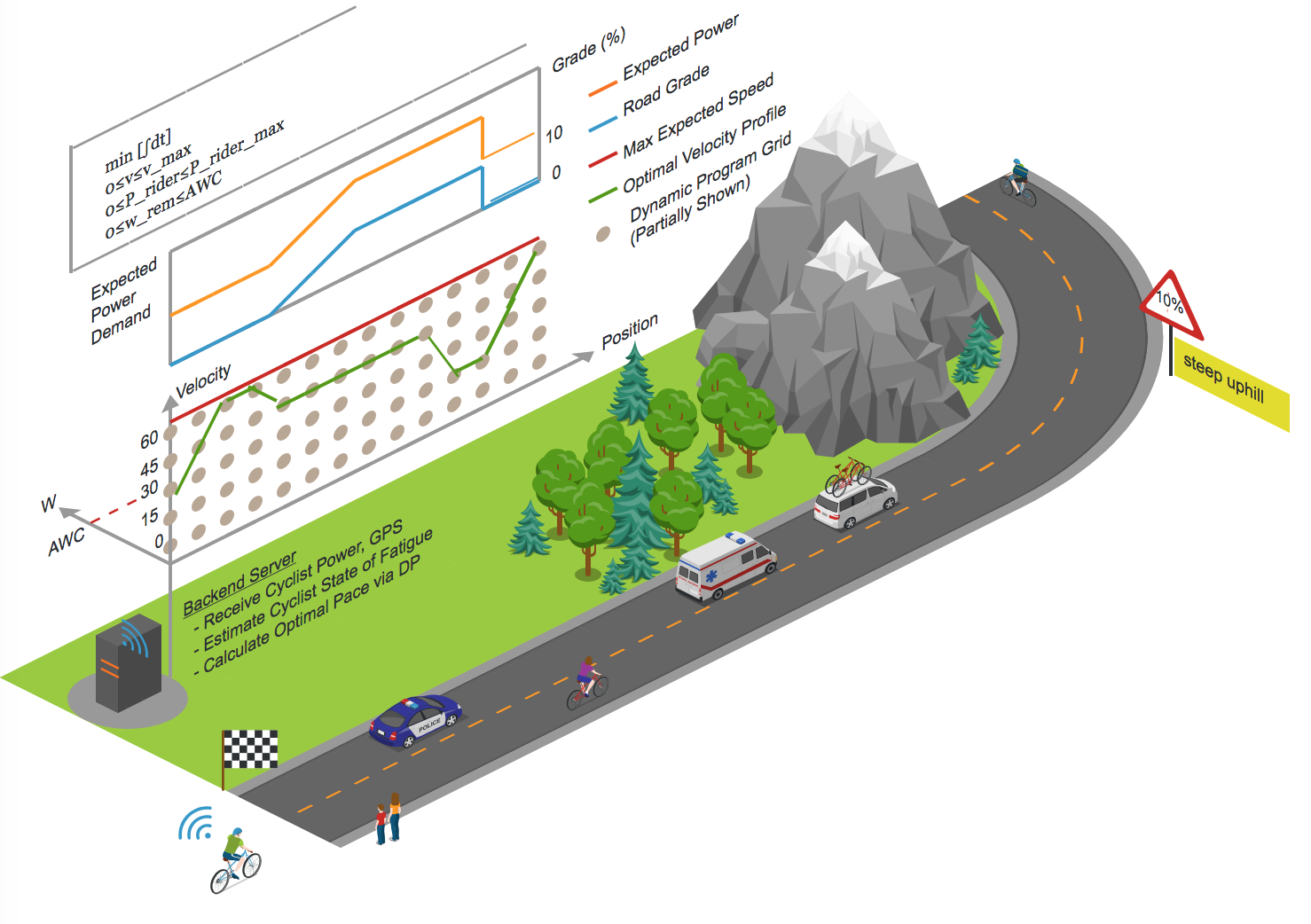}}
\caption{\footnotesize Illustration of an optimal feedback controller for pacing a cyclist during a time trial effort. Drawn by https://www.icograms.com.}
\label{fig_system}
\end{figure}

There are a few papers that addressed the optimal cycling problem. Authors in \cite{fayazi2013optimal} and \cite{wan2014optimal} suggest dynamic models for muscle's fatigue and recovery of force generation capacity. The models are used to solve an optimal control problem using dynamic programming to pace a cyclist during a time trial effort. However, the proposed models are not derived from experiments. In \cite{de2017individual}, the authors propose a muscle energy consumption model based on the critical power concept. Then, they formulate an optimal control problem for a 5km time trial on a flat road without considering energy recovery of muscle. The optimization problem is solved by applying Pontryagin's maximum principle.

In this paper, we propose the optimal strategy for a time trial cycling effort. In order to do so, we designed and executed a series of human subject experiments to develop muscle fatigue and recovery models that are explained in Section \ref{exp_protocol}. In Section \ref{modeling}, the experimental results and the models developed for one of our subjects is presented. Section \ref{optimization} explains the optimal control formulation which is solved using dynamic programming. The simulation results are presented in Section \ref{results}.

\section{Experimental Protocol}
\label{exp_protocol}
Fatigue has been commonly investigated in several research studies. There are several different definitions for fatigue according to those papers. For example, in \cite{Vollestad1997measurement}, fatigue is defined as any exercise-induced reduction in the maximal capacity to generate force or power. In another study \cite{edwards1981human} the author defines fatigue as the inability to maintain a required or desired force or power. Authors in \cite{sahlin1998energy} suggest that fatigue is determined from metabolic systems' depletion and signs of energy deficiency. Moreover, there is not a global mechanism for muscle fatigue. According to the findings in \cite{enoka2008muscle}, fatigue mechanism is specific to the task done by a human. Our focus in this research study is on cycling, so we will summarize a definition specifically for this task.

Although there are several scientific papers investigating muscle fatigue, there is not much on muscle recovery, especially while performing a task. Some of the studies such as \cite{bogdanis1995recovery} focus on recovery of metabolic system after finishing a physical task, and some other focus on the effect of personal diet on recovery \cite{saunders2004effects}. In \cite{shearman2016modeling} authors present a mathematical model for muscle anaerobic energy recovery. They suggest that energy recovery rate is not the same as energy expenditure, and that it is dependent on both recovery power level and duration.

In this study, we present the following definitions for muscle fatigue and recovery during cycling:

\begin{itemize}
\item \textit{Fatigue}: Expending energy from anaerobic metabolic systems by pedaling above critical power which results in a decrease of maximum power generation ability.
\item \textit{Recovery}: Recovering anaerobic metabolic systems energy by pedaling below critical power which results in an increase of maximum power generation ability. 
\end{itemize}

\begin{figure}
\centering
\setlength\fboxsep{0pt}
\setlength\fboxrule{0pt}
\setlength\belowcaptionskip{-20pt}
\fbox{\includegraphics[width=5 cm]{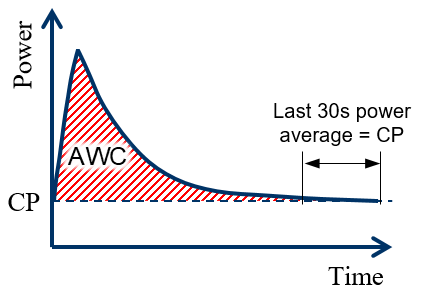}}
\caption{\footnotesize 3-minute-all-out test protocol. The average power at the last 30 seconds of the tests is considered to be $CP$, and the area between power plot and $CP$ is equivalent to $AWC$.}
\label{fig_3mao}
\end{figure}

According to the definitions above, we need to design a set of experiments by which we can develop mathematical models for expending and recovering energy from $AWC$ source, and also determine maximum power generation capacity of cyclists at any point during a ride. The tests are done in laboratory environments located in both Clemson University and Furman University campuses. Human Subjects are recruited from the cycling communities around Clemson and Greenville area. CompuTrainer from RacerMate is used to do the tests \cite{computrainer}. The trainer can be programmed to apply enough resistance to the rear wheel so the subjects apply the required power in our protocol. The experiment protocol includes thirteen lab visit days per subject.
On the first day, subject is instructed to perform a $\dot{V}_{O_2}$ ramp test. $\dot{V}_{O_2}$ is the volume of oxygen that a subject breathes in at any instance of time during the test. In this test power is increased incrementally until exhaustion. At a certain point during the test, there is a sudden change in the increase rate of $\dot{V}_{O_2}$ where a subject starts to burn energy from the anaerobic reservoir. This change point is called the Gas Exchange Threshold. We measure the power at this threshold to design our interval tests which are explained later on. The subjects do the 3-min-all-out test as in figure \ref{fig_3mao} twice so that we can better ensure consistency in our data.

\begin{figure*}
\centering
\setlength\fboxsep{0pt}
\setlength\fboxrule{0pt}
\setlength\belowcaptionskip{-20pt}
\fbox{\includegraphics[width= 11 cm] {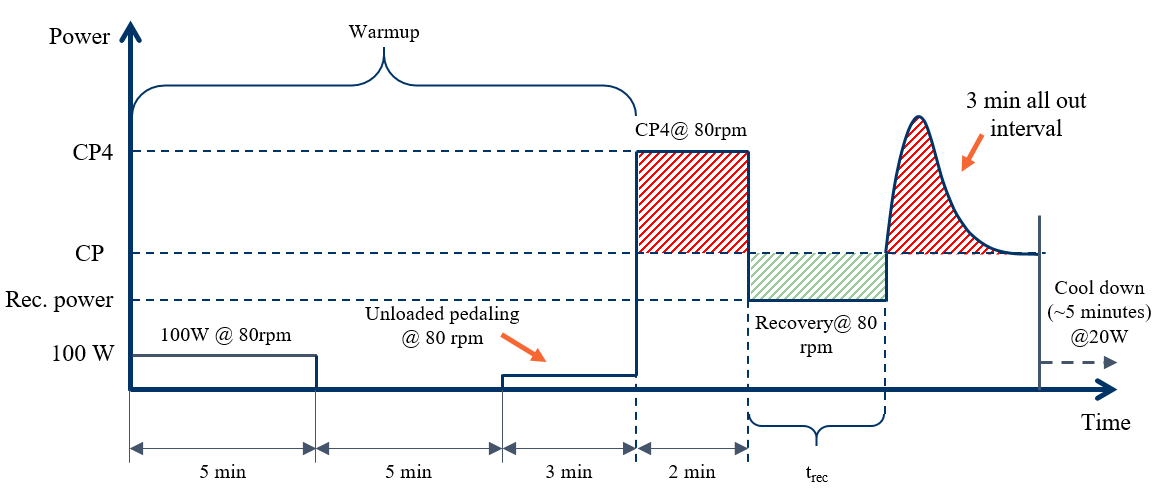}}
\caption{\footnotesize The power interval test protocol. After a warm-up period, the subject will pedal at $CP4$ for 2 minutes. Then he or she will pedal at three different recovery power levels for different time intervals to recover energy. After that, the subject will perform a 3-min-all-out to burn all the remaining energy from $AWC$.}
\label{fig_interval}
\end{figure*}

Also, there are 9 interval tests to understand recovery of $AWC$. Figure \ref{fig_interval} demonstrates this test. At the start of each test, subjects go through the warm-up protocol. Then the subjects are required to pedal for 2 minutes at $CP4$. $CP4$ is the power level the subject can maintain for 4 minutes without any drop in power. It is calculated by dividing $AWC$ by 240 seconds using equation (\ref{eq_basic-model}). By pedaling at $CP4$ for two minutes, the subjects are expected to lose half of their anaerobic energy. After that, the recovery interval begins. In order to model recovery, we perform the test at three different recovery power levels and three time intervals, which results in nine interval tests for each subject. Then, the subject will perform a 3-min-all-out test to burn all of the remaining energy. From these tests, one can determine how much power was recovered during each recovery interval by summing the areas above $CP$ in the first two minutes and 3-min-all-out, and then subtracting $AWC$ from it. Also, we can understand how power level and interval time affect the amount of recovered energy.

In addition to power data, we collect data for a parameter called $SmO_2$, which is a measure of muscle oxygenation. Muscle oxygenation plays an important role in muscle's ability to generate power. $SmO_2$ is at its maximum when cyclist is fresh and decreases as the cyclist gets fatigued. $SmO_2$ is the ratio between the concentration of oxygenated hemoglobins to that of the total hemoglobins in the blood flow of the local muscles. $SmO_2$ is measured using a non-invasive sensor called Moxy \cite{moxy}. This sensor uses Near-Infrared Spectroscopy method to measure muscle oxygenation in muscle tissue.

\section{Experiment results and models}
\label{modeling}
We have conducted experiments on a total of 9 subjects so far. However, because of the complexities in scheduling test days with subjects, we have been able to collect the most data from one of our subjects. Here we present the test results for our 9th subject, whom we refer to as Sub 9 for the rest of this paper. Recovery power levels from minimum to maximum are: 1) Smallest possible load on CompuTrainer, which is around 80 Watts, 2) 90\% of the power at gas exchange threshold ($0.9P_{GET}$), which is the power at which the cyclist starts to burn energy from anaerobic reservoir, and 3) Half way between $P_{GET}$ and $CP$. The fatigue and recovery intervals are done at the 80 rpm cadence.

We propose applicable mathematical models for Sub 9's energy consumption/recovery and maximum power generation ability. As discussed in Section \ref{introduction}, $AWC$ serves as a tank of energy while a cyclist is pedaling above critical power. When this tank is empty a cyclist can pedal only at $CP$. We define $W$ as the remaining energy from $AWC$. The rate of change of $W$ while expending energy is defined as the difference between rider's power and $CP$. One can assume that a similar equation can be used to calculate recovered energy. However, recovery happens at a slower rate than fatigue as shown in \cite{ferguson2010effect} and \cite{bickford2018modeling}. Therefore, the subject's actual recovered energy is less than the area in recovery portion of figure \ref{fig_interval}. Therefore, we can divide the actual recovered energy by the recovery time to obtain an adjusted recovery power:

	\begin{equation}
	P_{adj} \ =\ CP-\frac{W_{rec}}{T_{rec}}
	\label{eq_prec}
	\end{equation}

where $P_{adj}$ is the adjusted recovering power as determined by experimental data and $T_{rec}$ is the duration of recovery. We can rewrite equation(\ref{eq_prec}) for energy recovery along with the energy expenditure equation to form the switching energy model as below:


	\begin{equation}
    \frac{dW}{dt}=\left\{
    \begin{array}{@{} ccc @{}}
       -(P-CP) & P > CP & \text{(a)}\\
      \\
       -(P_{adj}-CP) & P < CP & \text{(b)}
    \end{array}\right.
  \label{eq_p_exp_rec}
	\end{equation}

where $P$ is the applied power to the bicycle by the rider. 

\begin{figure}
\setlength\fboxsep{0pt}
\setlength\fboxrule{0pt}
\fbox{\includegraphics[width=\columnwidth]{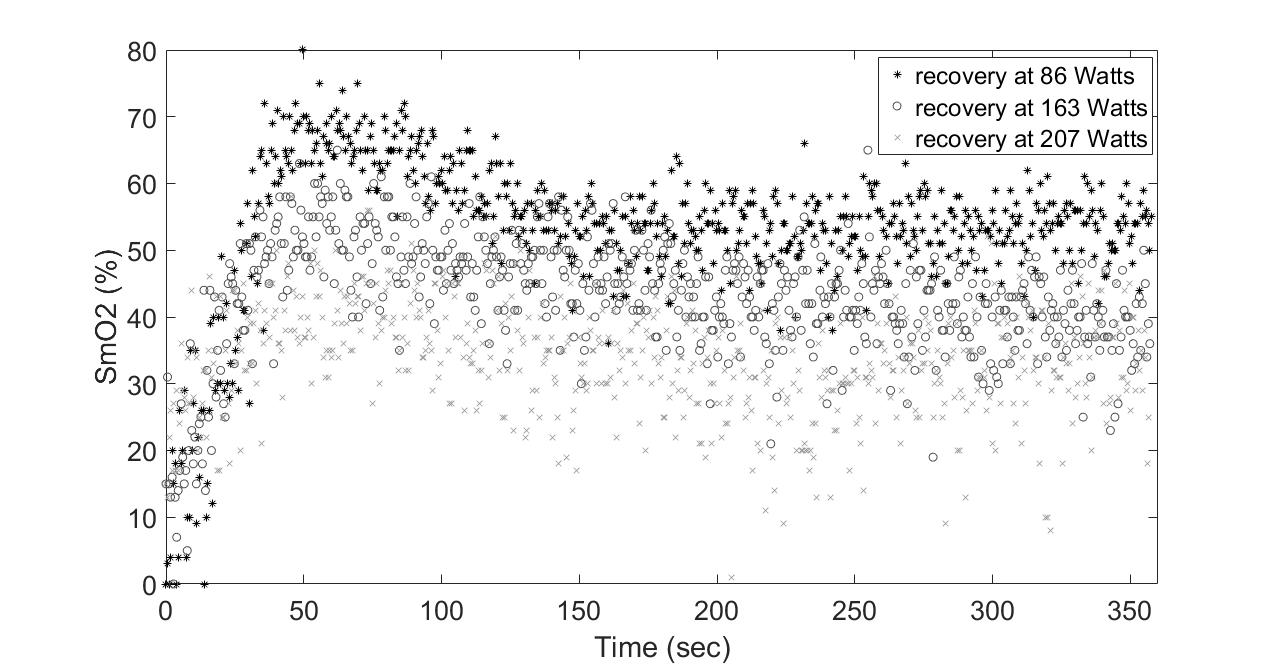}}
\caption{\footnotesize $SmO_2$ recovery during 6 min of recovery interval at three different power levels bellow sub 9's $CP$.}
\label{fig_smo2}
\end{figure}


As mentioned in Section \ref{exp_protocol}, one of the parameters measured during our tests is $SmO_2$. A high $SmO_2$ demonstrates the cyclist's high amount of remaining energy from $AWC$. Our data shows that, at each recovery power level, $SmO_2$ rises to a certain level regardless of the time period of recovery. The $SmO_2$ recovery is dependent on the power level as shown in figure \ref{fig_smo2}. Also, if the recovery model is based on interval time, it will not be a causal model. Therefore, we can take an average of adjusted powers (at 2 min, 6 min, 15 min recovery) at each recovery power level, and plot these adjusted powers vs. actual powers as in figure \ref{fig_rec-model} to develop our recovery model as:

\begin{figure}
\setlength\fboxsep{-5pt}
\setlength\fboxrule{-5pt}
\setlength\abovecaptionskip{-5pt}
\setlength\belowcaptionskip{-10pt}
\fbox{\includegraphics[width=\columnwidth]{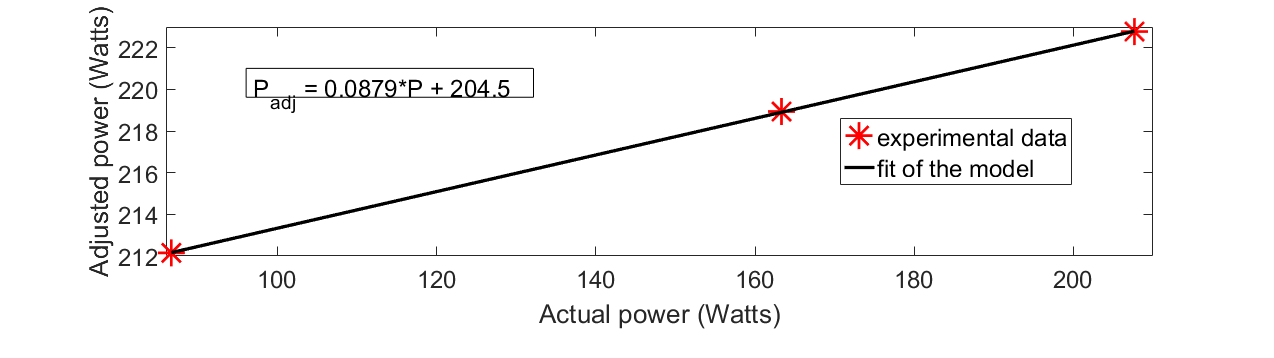}}
\caption{\footnotesize Adjusted recovery power vs. actual applied power by Sub 9 during the interval tests.}
\label{fig_rec-model}
\end{figure}
	\begin{equation}
	P_{adj} \ =\ aP + b
	\label{eq_p_adj}
	\end{equation}

where $a$ and $b$ are constants determined by fitting equation (\ref{eq_p_adj}) to experimental data, and are reported in figure \ref{fig_rec-model}, and $P$ is the recovery power of the cyclist.

Another parameter that is affected by cyclist's muscle fatigue and recovery is his or her maximum power generation ability at every point in time during cycling. In \cite{de2017individual}, authors considered a constant maximum power generation ability of 800 watts, which means the cyclist can always apply 800 watts and maintain it while he or she has energy, regardless of the history of power generation. This assumption is not valid as indicated by the 3-min-all-out test. As discussed in Section \ref{introduction}, a 3-min-all-out test requires subjects to apply maximum power at all times. However, their power decreases as they expend energy from their $AWC$. 

\begin{figure}
\setlength\fboxsep{-5pt}
\setlength\fboxrule{-5pt}
\setlength\belowcaptionskip{-10pt}
\fbox{\includegraphics[width=\columnwidth]{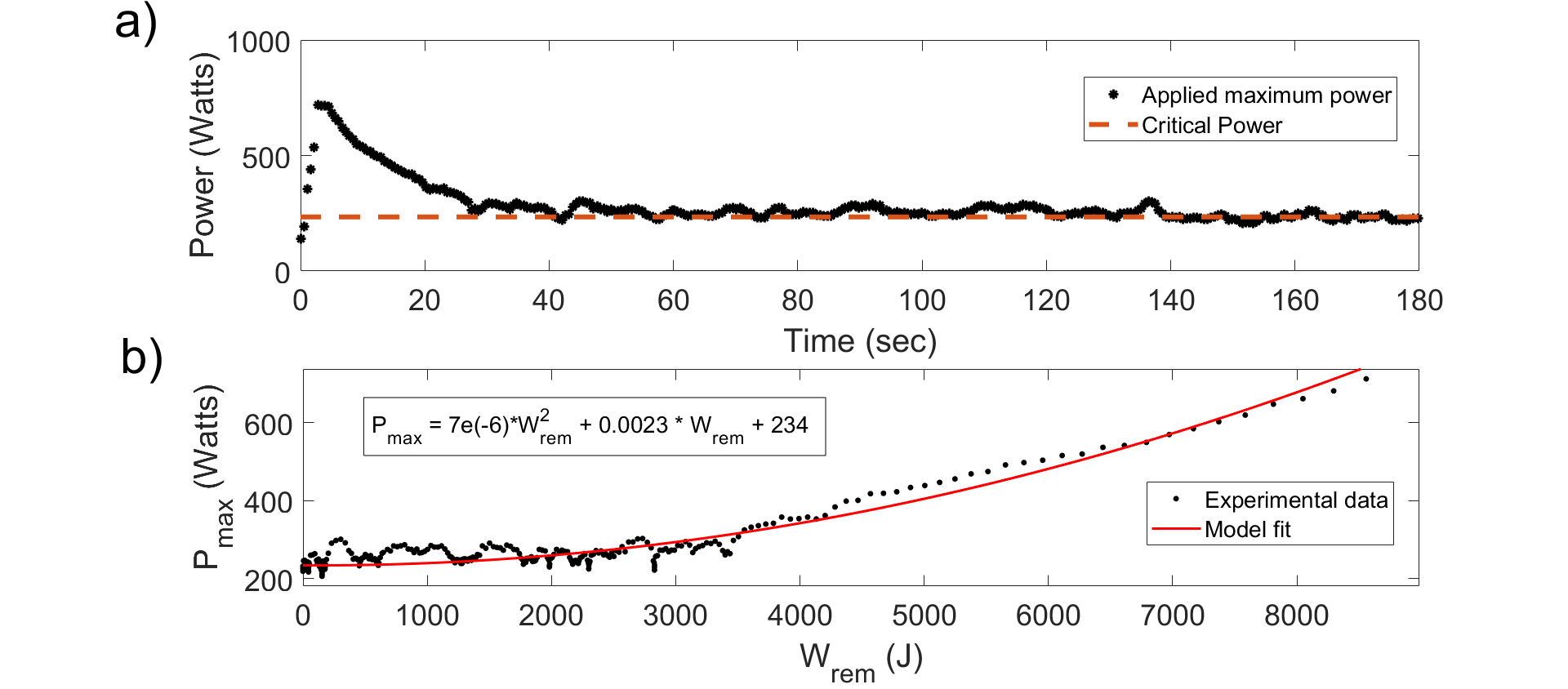}}
\caption{\footnotesize a) Sub 9's power data during the 3-min-all-out test. b) Maximum power model fit based on remaining energy during this test.}
\label{fig_3mao-experiment}

\end{figure}
A 3-min-all-out test as in figure \ref{fig_3mao-experiment}a can be used to construct a model of amount of $AWC$ remaining ($W$) vs. the maximum power that can be produced. By using the definitions presented in Section \ref{introduction}, Sub 9's $CP$ and $AWC$ are calculated to be 234 Watts and 9758 J, respectively. We can then calculate the remaining anaerobic energy at any point in time during the test by calculating the remaining  area between power plot and $CP$. Then we plot the remaining energy at every time vs. maximum power applied at that particular time as in figure \ref{fig_3mao-experiment}b. For creating this plot, the data points before the maximum power are removed, because the subject is trying to overcome wheel and aerobic inertia to reach the maximum power. We then fit the maximum power model to a quadratic expression:

	\begin{equation}
	P_{max}(t) \ = \ a_1W^2(t) + a_2W(t) + CP
	\label{eq_max-p}
	\end{equation}

where $W$ is the remaining energy form the cyclist's available work capacity, and $a_1$ and $a_2$ are model constants calculated from experimental data in figure \ref{fig_3mao-experiment}b.

\section{Optimal Control Formulation} 
\label{optimization}

In this section, the problem formulation is discussed. According to Newton's Second Law, the dynamics of bicycle motion can be formulated as below:

	\begin{equation}
	\begin{split}
	&F(t) \ = \ m_t \frac{dv(t)}{dt} + m_t g (sin(\theta) + \mu cos(\theta)) \\
	& + 0.5 C_d \rho A v(t)^2 
	\label{eq_newton_f}
	\end{split}
	\end{equation}

Where $F(t)$ is the bicycle's driving force at its wheels, $m_t$ is the total mass of the bicycle and rider, $C_d$ is the drag coefficient, $A$ is the frontal area, $\rho$ is the density of air which is assumed to be constant and independent of the elevation, $\theta$ is the road slope angle which is positive for uphill and negative for downhill, and $\mu$ is the rolling resistance of the road.

Since our muscle models are based on rider power, equation (\ref{eq_newton_f}) can be reformulated by driving power $P(t)=F(t)v(t)$ as in \cite{de2017individual}:

	\begin{equation}
	\begin{split}
	&P(t) \ = \ \bigg(m_t \frac{dv(t)}{dt} + m_t g (sin(\theta) + \mu cos(\theta)) \\
	& + 0.5 C_d \rho A v(t)^2 \bigg) v(t)
	\label{eq_newton_p}
	\end{split}
	\end{equation}
	
Assuming 100\% efficiency for bicycle powertrain, we can assume $P(t)$ is the cyclist's power on the pedals. The advantage of using equation (\ref{eq_newton_p}) is that gear selection is not a factor in our formulation, which otherwise makes the optimization more complex.

An optimal control problem is formulated to calculate bicycle velocity over time, such that the rider finishes a course in minimum time. Because cycling path is known, while final time is dynamic, the objective function is formulated and reparameterized as:

	\begin{equation}
	\text{min} \ {J = \int_{0}^{t_f} dt = \int_{0}^{x_f} \frac{1}{v(t)}dx}
	\label{eq_t}
	\end{equation}

The optimization constraints include limits on power, velocity and energy as below:

	\begin{equation}
    \left\{
    \begin{array}{@{} l c @{}}
      \text{power limit: } & 0 \leq P_{rider}(x) \leq P_{rider,\,max}(x) \\
      \text{remaining energy limit: } & 0 \leq W(x) \leq AWC \\
      \text{velocity range: } & 0 \leq v(x) \leq v_{max}
    \end{array}\right.
  \label{eq_constraint}
	\end{equation}

where $W$ is the remaining energy of the cyclist which was previously defined as $W_{rem}$.

Now the optimization formulation has to be discretized to be solved by dynamic programming (DP). Velocity and remaining energy are the two states of the system, and distance is the independent variable that we base the discretization on. We can combine equations (\ref{eq_p_exp_rec}) and (\ref{eq_p_adj}) to write the fatigue and recovery models in discrete time as below:

	\begin{equation}
    \left\{
    \begin{array}{@{} l c @{}}
      P_{rider_{,\,i}}  = -(W_{i+1} - W_i)\frac{v_i+v_{i+1}}{2\Delta X} + CP& W_{i+1} < W_i\\
      \\
      P_{rider_{,\,i}}  = -\frac{1}{a}((W_{i+1} - W_i)\frac{v_i+v_{i+1}}{2\Delta X} - b + CP)& W_{i+1} > W_i
    \end{array}\right.
  \label{eq_p_rider}
	\end{equation}
	
Where $P_{rider}$ is the cyclist's power output, $dX$ is the distance between two consecutive states, and $i$ and $i+1$ are current state and next state, respectively. Also, $dt$ is substituted with $\frac{2\Delta X}{v_i+v_{i+1}}$ since our independent variable is chosen to be distance $X$. In addition to rider's power, we can calculate the required power from equation (\ref{eq_newton_p}) as below:

\begin{equation}
	\begin{split}
	&P_{req_{,\,i}} \ = \ \bigg(m_t \frac{v_{i+1}^2 - v_i^2}{2\Delta X} + m_t g (sin(\theta_i) + \mu cos(\theta_i)) \\
	& + 0.5 C_d \rho A \frac{(v_{i+1} + v_i)^2}{2} \bigg) \frac{v_{i+1} + v_i}{2}
	\label{eq_newton_p_dis}
	\end{split}
\end{equation}

The rider's power from equation (\ref{eq_p_rider}) should match the required power from equation (\ref{eq_newton_p_dis}) to get a feasible result. We add a constraint in our DP code to handle this requirement. 

In summary, the state-space model can be written as below:

	\begin{equation}
	z(i+1) \ =\ f(z(i),u(i))
	\label{eq_states_DP}
	\end{equation}
	
where $z=\begin{bmatrix}
v\\
W
\end{bmatrix}$ are the two states of the system while $u=P_{rider}$ is the input to it. Also, we have an equality constraint on input $P_{rider}=P_{req}$ along with the in-equality constraint in equation (\ref{eq_constraint}).
\section{Optimization Results}
\label{results}
\subsection{Baseline Test}
\begin{figure}
\setlength\fboxsep{0pt}
\setlength\fboxrule{0pt}
\setlength\belowcaptionskip{-20pt}
\fbox{\includegraphics[width=\columnwidth]{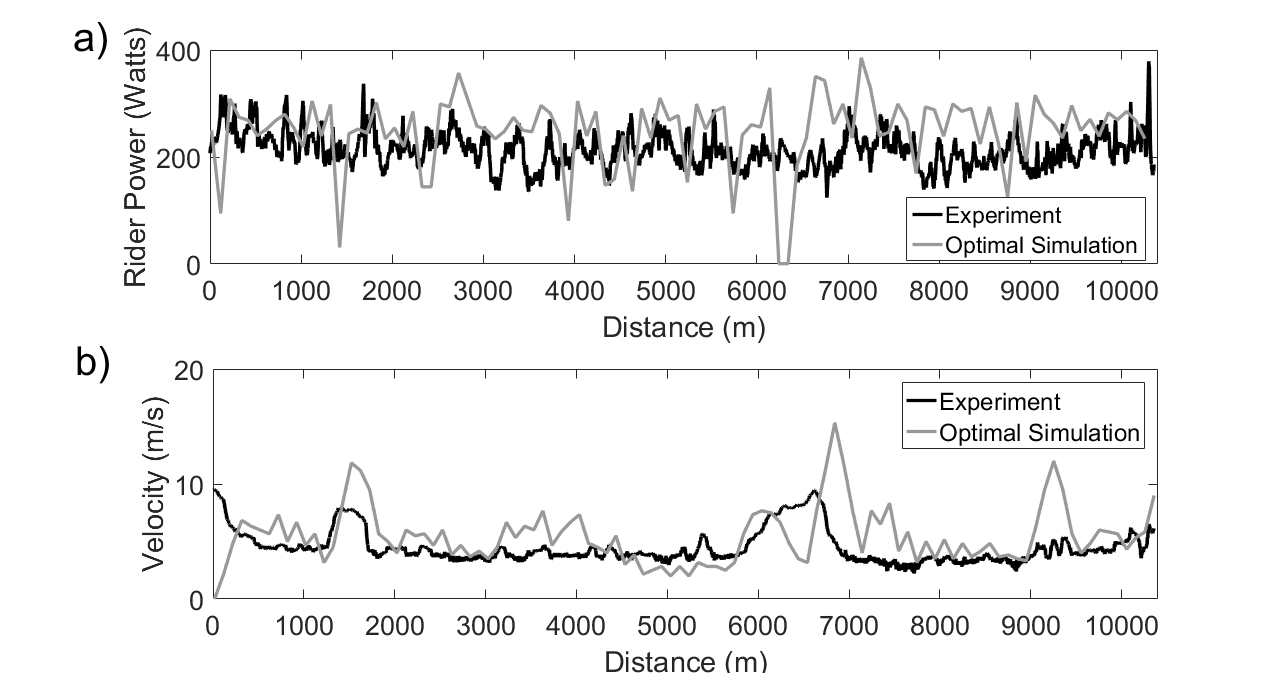}}
\caption{\footnotesize Comparison between subject's power and velocity during the baseline test and optimal simulation.}
\label{fig_optimal-vs-baseline}
\end{figure}

Now that we have our fatigue and recovery models along with optimal control formulation we can use DP to solve the optimal control problem. Sub 9 chose a course that he/she rides through frequently. The chosen course was the cycling track in Ceaser's Head state park, Greenville SC. The course is 10.3 km long, and includes two long uphills with a short flat and slightly downhill section in between. The course was simulated on CompuTrainer using PerfPro studio software by providing elevation data to it. As a baseline for comparison, Sub 9 was asked to ride on the CompuTrainer and finish the simulated course in shortest time possible using his/her own strategy, and was verbally encouraged to ride as quickly as possible. Sub 9 finished the test in 41 minutes and 17 seconds. Since the subject was already familiar with the course, they verbally verified that the CompuTrainer simulated course was very close to the real one. The subject's power and velocity data during the test can be seen in figure \ref{fig_optimal-vs-baseline}. 

\subsection{Optimal Control Simulation}
It should be noted that in the laboratory environment, drag force does not exist. Therefore, we should remove the drag term in equation (\ref{eq_newton_p_dis}) in our dynamic programming formulation.

One of the challenges in writing the DP code is to the adjusted recovery power from equation \ref{eq_p_adj}. According to the values for $a$ and $b$ from figure \ref{fig_rec-model}, for a specific interval of power output above $CP$, the adjusted power $P_{adj}$ can have a value below $CP$. In this condition, the subject will recover energy while pedaling above $CP$, which is not reasonable. Therefore, we added an upper limit constraint on applied power at $CP$ when the subject is recovering energy.

\begin{figure}
\setlength\fboxsep{0pt}
\setlength\fboxrule{0pt}
\setlength\abovecaptionskip{-10pt}
\setlength\belowcaptionskip{-20pt}
\fbox{\includegraphics[width=\columnwidth]{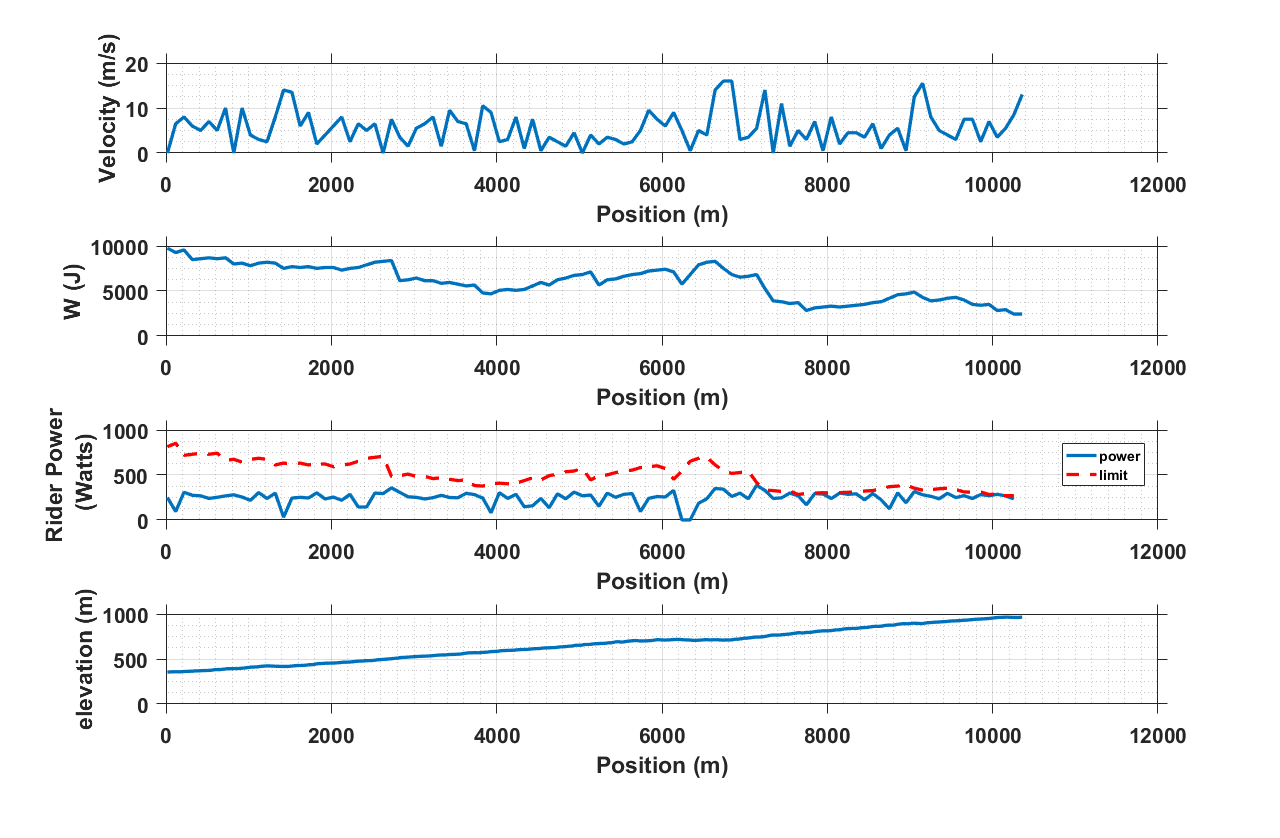}}
\caption{\footnotesize Optimization simulation results for a 10.3 km road in Greenville SC.}
\label{fig_optimal}
\end{figure}

By considering maximum velocity of $16\,m/s$, a grid of 32 velocity nodes by 100 remaining energy nodes along with 100 m distance intervals is constructed. The results of the optimal control simulation can be seen in figure \ref{fig_optimal}. According to the results, the subject would be able to finish the same course in 38 minutes and 15 seconds. The power and velocity of the subject during the baseline test and the optimal simulation are compared in figure \ref{fig_optimal-vs-baseline}. We can observe that the subject's own strategy does not involve recovery at very low power levels. However, in the optimal solution, by recovering more energy during recovery periods, the subject is able to pedal at higher powers for longer time periods, especially during the last kilometer of the course. 

Unlike the computation power used by authors in \cite{fayazi2013optimal}, the grid size used in this paper helped us to run the DP code on a regular desktop computer which has an Intel(R) Core(TM) i5-4460 CPU at 3.2 GHz, and 12 GB of RAM. The entire simulation takes less than one minute to run for the 10.3 km course. This low computational load will enable us to develop a real-time software-in-the-loop controller.

\subsection{Potential for Implementing the Controller}
Our ultimate goal for this project is to implement our optimal feedback controller in an experiment to reduce the cyclist's travel time on a course. For implementation, the DP is to be rewritten as a real-time controller. The goal is to develop a software-in-the-loop test, in which a subject is provided with the optimal power to hold at each position during the test. In order to develop the software, the backward portion of the DP is solved and its matrices are saved. These matrices can be used as a look-up table for performing the forward DP in real-time. The states of the system (velocity and remaining energy) are estimated using real-time data logged from the CompuTrainer. Then, using the backward DP matrices we can calculate the optimal power to hold for the cyclist. Since there is variation in power output, the rider may not be able to apply the exact optimal power during each distance interval. Therefore, the real-time power data is recorded and fed back to the software. The average applied power during each distance interval is then calculated and used for estimating the updated states of the system. The estimation enables us to update the optimal solution in real-time. 

In the near future we plan to perform the optimal test with the subjects we are currently recruiting in order to assess the applicability of our optimal solution for time trial efforts.  
%
%
\section{Conclusions}
In this paper we address a scientific challenge in modeling and optimizing a cycling time trial effort. Utilizing the concepts of \textit{Anaerobic Work Capacity} and \textit{Critical Power}, we developed an experimental protocol to mathematically model a cyclist's fatigue and recovery. The models are based on the applied power by the cyclist, and the corresponding decrease or increase in his or her remaining energy. The fatigue and recovery models are then used to formulate an optimal control problem for a time trial on a course in Greenville SC. As a baseline for comparison, the subject was asked to ride on the course simulated on the CompuTrainer in a laboratory environment, and use his/her own strategy to finish it in shortest time possible. Then, the optimal control problem was solved using dynamic programming. The optimization results show 3 minutes reduction of travel time compared to the baseline test. In optimal simulation, the subject can benefit from recovery intervals at low powers to regain energy which can be expended to generate higher power levels, especially towards the end of the road. In addition to the simulation, a software-in-the-loop test is developed to implement a real-time controller which uses power data as feedback. The rider will be provided with optimal power to hold during each distance interval, in our future work.

%

\section*{ACKNOWLEDGMENT}
The authors want to thank Furman University students Mason Coppi, Frank Lara, Lee Shearer, Nicholas Hayden, Jake Ogden, and Brendan Rhim for their contribution in data collection. 

%
\bibliographystyle{ieeetr}
\bibliography{root}
\end{document}